\documentclass{scrartcl}

\usepackage[anythingbreaks]{breakurl}
\usepackage{jmbarrett}
\usepackage{tikz-cd}
\tikzcdset{row sep/normal=2.7em,column sep/normal=3.6em}

\nocite{*}

\renewcommand{\C}{\mathbf{C}}
\newcommand{\intl}{\mathbf{0}}
\newcommand{\tl}{\mathbf{1}}
\newcommand{\unq}{\textbf{!}}
\newcommand{\mon}{\hookrightarrow}
\newcommand{\SC}{\Omega}
\newcommand{\true}{\textsc{t}}
\newcommand{\false}{\textsc{f}}
\newcommand{\ETCS}{\mathsf{ETCS}}
\newcommand{\BZC}{\mathsf{BZC}}
\renewcommand{\R}{\mathsf{R}}
\newcommand{\ch}{\textsc{c}}
\newcommand*\rcolon{%
	\nobreak
	\mskip6mu plus1mu
	\mathpunct{}%
	\nonscript
	\mkern-\thinmuskip
	{:}%
	\mskip2mu
	\relax
}

\DeclareMathOperator{\Sub}{Sub}
\DeclareMathOperator{\Hom}{Hom}
\DeclareMathOperator{\Pp}{\mathcal{P}\!}
\DeclareMathOperator{\ev}{ev}
\DeclareMathOperator{\im}{im}
\DeclareMathOperator{\coim}{coim}

\title{Elementary topoi}
\author{Jordan Mitchell Barrett}
\date{November 13, 2020}

\begin{document}

\maketitle

\section{Introduction}

As the prototypical category, $\Set$ has many properties which make it special amongst categories. From the point of view of mathematical logic, one such property is that $\Set$ has enough structure to ``properly'' formalise logic. However, we could ask what it might mean to formalise logic in another category $\C$. The notion of an (\textit{elementary}) \textit{topos} distills the essential features of $\Set$ which allow us to do this.

Throughout this report, a boldface $\C$ will denote an arbitrary category (and later, a topos). Objects of $\C$ are denoted by uppercase letters $A, B, C, \ldots$, and morphisms by lowercase $f,g,h,\ldots$. We assume familiarity with basic category theory, but let us recall some of the important categorical notions we will need.

\begin{definition}\label{defn:monic-epi}
	An arrow $f\colon B \to C$ in $\C$ is \textit{monic} (denoted $f\colon B \mon C$) if for all arrows $g, h\colon A \to B$, $fg = fh$ implies $g=h$. Dually, $f$ is \textit{epi} if for all arrows $g, h\colon C \to D$, $gf = hf$ implies $g=h$.
\end{definition}

\begin{definition}
	A \textit{pullback} is a limit over a diagram of shape \begin{tikzcd}[cramped,sep=small] \bullet \arrow[r] & \bullet & \bullet \arrow[l] \end{tikzcd}. More concretely, the pullback of two maps $f\colon B \to D \leftarrow C \rcolon g$ is a pair $h\colon B \leftarrow A \to C \rcolon k$ making the square commute:\vspace*{-5mm}%
	\begin{center}%
	\begin{tikzcd}%
		A \arrow[d,"h"] \arrow[r,"k"] & C \arrow[d,"g"] \\
		B \arrow[r,"f"] & D
	\end{tikzcd}
	\end{center}
	and which is \textit{universal}---any other such pair $h'\colon B \leftarrow A' \to C \rcolon k'$ factors through $h,k$ via a unique $t\colon A' \to A$:\vspace*{-3mm}
	\begin{center}
	\begin{tikzcd}[sep=small]
		{\color{darkgray} A'} \arrow[ddddr,bend right,"h'"',shorten >=-1mm,darkgray] \arrow[rrrrd,bend left,"k'",shorten=-1mm,darkgray] \arrow[dr,"t",dashed,shorten=-1.5mm,darkgray] \\
		& A \arrow[ddd,"h"] \arrow[rrr,"k"] & & & C \arrow[ddd,"g"] \\
		\\ \\
		& B \arrow[rrr,"f"] & & & D
	\end{tikzcd}
	\end{center}
\end{definition}

Calling it \textit{the} pullback is justified, since pullbacks are unique up to unique isomorphism. A useful fact is that if $f$ is monic, its pullback $k$ is monic too \cite[Ex.\ III.4.5]{maclaneCategoriesWorkingMathematician1998}.

\begin{definition}\label{defn:terminal}
	An object $\intl \in \C$ is \textit{initial} if for every object $A \in \C$, there is a unique arrow $\unq\colon \intl \to A$. Dually, $\tl \in \C$ is \textit{terminal} if there is a unique arrow $\unq\colon A \to \tl$ from any object $A$.
\end{definition}

For example, the empty set $\varnothing$ is initial in $\Set$, while any singleton set $\{ * \}$ is terminal. The trivial group $\{ 0 \}$ is both initial and terminal in $\Grp$ and $\Ab$. Definitions \ref{defn:monic-epi} and \ref{defn:terminal} imply any map $\tl \to B$ is monic.

\section{Subobject classifiers}

Henceforth, we will assume that $\C$ has all finite limits and colimits. This implies that $\C$ has an initial object $\intl$ and terminal object $\tl$, since they are the colimit and limit, respectively, of the empty diagram.  For $X \in \C$, $X^n$ denotes the $n$-fold categorical product.

The object $\mathbf{2} = \{ 0, 1 \}$ in $\Set$ plays a very special role, in that for any set $X$, there is a correspondence between subsets $Y \subseteq X$ and characteristic functions $\chi_Y\colon X \to \mathbf{2}$. In more general categories, this idea is captured by the notion of a \textit{subobject classifier}.

\begin{definition}\label{defn:subobj-cl}
	A \textit{subobject classifier} is an object $\SC \in \C$, along with a morphism $\true\colon \tl \to \SC$, such that for every monic $f\colon A \mon B$ in $\C$, there is a unique map $\chi_f\colon B \to \SC$, called the \textit{character} of $f$, making the following a pullback:
	\begin{center}
	\begin{tikzcd}
		A \arrow[r,"\unq"] \arrow[d,"f"',hook] & \tl \arrow[d,"\true"] \\
		B \arrow[r,"\chi_f",dashed] & \SC
	\end{tikzcd}
	\end{center}
\end{definition}

This definition can be recast in a more ``categorical'' way, using the intuition that $\SC$ determines a natural isomorphism between $\Hom(B,\SC)$ and the ``subobjects'' of $B$.

\begin{remark}
	For $B \in \C$, its \textit{subobjects} $\Sub(B)$ are the monics $f\colon A \mon B$ for some $A$, modulo isomorphisms of the domain, i.e.\ $f\colon A \mon B$ and $f'\colon A' \mon B$ are equivalent if there is invertible $g\colon A \to A'$ making the following commute:	\vspace*{-2mm}
	\begin{center}
	\begin{tikzcd}[row sep=small]
		A \arrow[rd,"f",hook] \arrow[dd,shift right,dashed,"g"'] & \\
		& B \\
		A' \arrow[ur,"f'"',hook] \arrow[uu,shift right,dashed,"g^{-1}"'] &
	\end{tikzcd}
	\end{center}

Define the (contravariant) \textit{subobject functor} $\Sub\colon \C \to \Set$ by $B \mapsto \Sub(B)$, and $f\colon B' \to B$ maps to $Sf\colon \Sub(B) \to \Sub(B')$, the pullback along $f$:\vspace*{-2mm}
	\begin{center}
	\begin{tikzcd}
		A \arrow[r,"g",hook] & B \\
		A' \arrow[r,"Sf(g)"',hook] \arrow[u] & B' \arrow[u,"f"']
	\end{tikzcd}
	\end{center}
	This is well-defined on equivalence classes. Then, a subobject classifier can equivalently be defined as a representation $(\SC, \varphi)$ of $\Sub$, where $\varphi\colon \Sub(-) \cong \Hom(-,\SC)$ is a natural isomorphism. The unique map $\id_\tl\colon \tl \to \tl$ is monic, so a subobject of $\tl$, whence the truth map can be defined $\true = \varphi_\tl(\id_\tl)$ \cite[\S2]{gokhmanSubobjectClassifierCategory1997,leinsterInformalIntroductionTopos2011}.
\end{remark}

As expected, the object $\textbf{2}$ is indeed a subobject classifier for $\Set$. Let's see a more interesting example.

\begin{example}\label{exm:slice-sc}
	For an object $X \in \C$, the \textit{slice category} $\C/X$ has maps $f\colon A \to X$ as objects, where morphisms $(f\colon A) \to (f'\colon A')$ are maps $g\colon A \to A'$ making the following commute:\vspace{-5mm}
	\begin{center}
	\begin{tikzcd}[row sep=small]
		A \arrow[rd,"f"] \arrow[dd,"g"'] & \\
		& X \\
		A' \arrow[ur,"f'"']
	\end{tikzcd}
	\end{center}
	$\C/X$ has a terminal object $\id_X\colon X \to X$. Now, if $\C$ has subobject classifier $(\SC,\true)$, then $\C/X$ has subobject classifier $\pi_2\colon \SC\!\times\!X \to X$, the canonical projection from the product, with truth map $\bar{\true} = \true \!\times\! \id_X\colon (\id\colon X) \to (\pi_2\colon \SC\!\times\!X)$. To see this, a map $f\colon (a\colon A) \to (b\colon B)$ is monic iff $f\colon A \to B$ is. In $\C$, we get $\chi_f\colon B \to \SC$; then, the composite\vspace*{-4mm}
	\begin{center}
	\begin{tikzcd}
		B \arrow[r,"\Delta_B"] & B^2 \arrow[r,"\chi_f \times b"] & \SC\!\times\!X
	\end{tikzcd}
	\end{center}
	is the character of $f$ in $\C/X$, where $\Delta_B\colon B \to B^2$ is the diagonal.
\end{example}

$\SC$ should be interpreted as an object containing all possible ``truth values'' for the category $\C$---we will develop this interpretation in Section \ref{sec:mod}.



\section{Topoi}

We are nearly ready to give the definition of a topos. We require the extra structure of \textit{cartesian closure}, meaning for any two objects $A, B \in \C$, there is an object $B^A$ which acts like the set of functions $A 
\to B$. The formal definition follows.

\begin{definition}\label{defn:expon}
	The \textit{exponential} of two objects $A, B \in \C$ is an object $B^A$ and a morphism $\ev\colon B^A \!\times\! A \to B$, such that for any morphism $g\colon C \!\times\! A \to B$, there is unique $\bar{g}\colon C \to B^A$ making the following commute:
	\begin{center}
	\begin{tikzcd}
		C \arrow[d,dashed,"\bar{g}"'] & C \!\times\! A \arrow[d,dashed,"\bar{g} \times \id_A"'] \arrow[rd,"g"] \\
		B^A & B^A \!\times\! A \arrow[r,"\ev"'] & B
	\end{tikzcd}
	\end{center}
\end{definition}

The prototypical example is in $\Set$: $B^A = \{ f\colon A \to B \}$, with $\ev(f,a) = f(a)$. As with Definition \ref{defn:subobj-cl}, this has a more categorical interpretation---a natural isomorphism $\Hom(C \!\times\! A,B) \cong \Hom(C,B^A)$. We could alternatively define $B^A$ as a representation of the functor $\Hom(- \!\times\! A,B)\colon \C \to \Set$; $\ev$ is recovered by mapping $\id_{B^A}$ through the associated natural isomorphism (for $C = B^A$). Furthermore, if $A$ is such that $B^A$ exists for all $B$, then $(-)^A\colon \C \to \C$ is a right adjoint to the product functor $(-) \!\times\! A\colon \C \to \C$.

\begin{definition}
	A \textit{cartesian closed category} (ccc) is one with all finite products and exponentials. A \textit{ccc$^+$} is a ccc with all finite limits and colimits.
\end{definition}
  
\begin{definition}
	An (elementary) \textit{topos} (plural \textit{topoi} or \textit{toposes}) is a ccc$^+$ with a subobject classifier.
\end{definition}

We can verify $\Set$ is a ccc$^+$, hence a topos. Furthermore, we saw in Example \ref{exm:slice-sc} that whenever $\C$ has a subobject classifier, so does $\C/X$ for any object $X \in \C$. Even better, whenever $\C$ is a topos, so is $\C/X$. This fact is sometimes called the \textit{fundamental theorem of topoi}, and apart from showing that $\C/X$ has exponentials, the proof is straightforward \cite[\S IV.7]{maclaneSheavesGeometryLogic1992}.

There are many more examples of topoi; a particularly important one is the functor category $\Set^{\C^\mathrm{op}}$ for \textit{any} category $\C$, whose objects are called \textit{presheaves on $\C$}. Topoi can also be characterised by the existence of objects $\Pp X$ which act like the ``power set'' of an object $X$.

\begin{definition}\label{defn:power-obj}
	A \textit{power object} of an object $X \in \C$ consists of objects $\Pp X$, $K \in \mathbf{C}$ and a monic $\ni\colon K \mon \Pp X \!\times\! X$ such that for every monic $r\colon A \mon B \!\times\! X$, there is a unique morphism $\bar{\chi}_r\colon B \to \Pp X$ such that the following is a pullback:
	\begin{center}
	\begin{tikzcd}
		A \arrow[rr] \arrow[d,"r"',hook] & & K \arrow[d,"\ni",hook] \\
		B \!\times\! X \arrow[rr,"\bar{\chi}_r \times \id_X",dashed] & & \Pp X \!\times\! X
	\end{tikzcd}
	\end{center}
\end{definition}

As always, the motivating example is in $\Set$, where $\Pp X$ is the actual power set of $X$, and $\ni$ the usual membership relation. We have written $\ni$ backwards to match the definition of an exponential; we now expound the connection.

\begin{theorem}\label{thm:topos-power-obj}
	A ccc$^+$ $\C$ has a subobject classifier if and only if $\C$ has all power objects.
\end{theorem}

\begin{proof}
	\begin{itemize}
		\iffb $\SC \defeq \Pp\tl$ is the required subobject classifier. Taking $X=\tl$ in Definition \ref{defn:power-obj}, and mapping through the natural isomorphism $\tl \!\times\! B \cong B$ gives Definition \ref{defn:subobj-cl}, since it is necessarily true that $K \cong \tl$ \cite[Fact 2.4]{leinsterInformalIntroductionTopos2011}.
		
		\ifff $\Pp X \defeq \SC^X$ is the required power object. Pulling back $\true\colon \tl \to \SC$ along $\ev\colon \SC^X \!\times\! X \to \SC$ gives $\ni\colon K \mon \SC^X \!\times\! X$, which is monic since $\true$ is. Now, for any monic $r\colon A \mon B \!\times\! X$, Definition \ref{defn:subobj-cl} gives $\chi_r\colon B \!\times\! X \to \SC$, whence Definition \ref{defn:expon} gives the required $\bar{\chi}_r\colon B \to \SC^X$.\qedhere
	\end{itemize}
\end{proof}

\section{The internal logic of a topos}
\label{sec:mod}

For $A \in \C$, the definition of the product $A^2 = A \!\times\! A$ gives the \textit{diagonal} $\Delta\colon A \to A^2$ as the unique map making the following commute:
\begin{center}
\begin{tikzcd}
	A & A^2 \arrow[l,"\pi_1"'] \arrow[r,"\pi_2"] & A \\
	& A \arrow[ul,"\id_A"] \arrow[ur,"\id_A"'] \arrow[u,dashed,"\Delta"]
\end{tikzcd}
\end{center}

This extends to $n$-fold products in the obvious way. More generally, suppose $\sigma = (m_1,\ldots,m_n)$ is an $n$-tuple of positive integers, with $d = \max\{m_i: i \leq n\}$. The \textit{diagonal of signature $\sigma$} is the unique map $\Delta_\sigma\colon A^d \to A^n$ making the following commute for each $i \leq n$:\vspace*{-5mm}
\begin{center}
\begin{tikzcd}
	A^n \arrow[r,"\pi_i"] & A \\
	A^d \arrow[ur,"\pi_{m_i}"'] \arrow[u,dashed,"\Delta_\sigma"]
\end{tikzcd}
\end{center}
The regular diagonals $\Delta\colon A \to A^n$ have signature $(1,1,\ldots,1)$.

Henceforth, let $\C$ be a topos. We introduced topoi on the promise that they would allow a good formalisation of logic. In this section, we make good on that promise, by showing how first-order logic can be internalised in a topos. Using exponentials, topoi can also formalise higher-order logic---for simplicity, we will not explore this here.

We assume familiarity with basic first-order model theory, as in \cite{markerModelTheoryIntroduction2002, changModelTheory1990}. In particular, for $\Ll$ a language, $\Ll$-terms and $\Ll$-formulae have their usual definitions.

\begin{definition}
	An \textit{$\Ll$-structure $\M$ in $\C$} consists of the following data:
	\begin{enumerate}
		\item An object $M \in \C$ called the \textit{support} or \textit{domain};
		
		\item For every $n$-ary function symbol $f \in \Ll$, a morphism $f^\M\colon M^n \to M$;
		
		\item For every $n$-ary relation symbol $R \in \Ll$, a morphism $R^\M\colon M^n \to \SC$;
		
		\item For every constant symbol $c \in \Ll$, a morphism $c^\M\colon \tl \to M$.
	\end{enumerate}
\end{definition}

A well-known example is a \textit{group object} in a category $\C$, which is a structure in the language $\{ *, (\cdot)^{-1}, e \}$ satisfying certain axioms.

Now, we define interpretations of terms and formulae in $\M$. Essentially, a term with $n$ free variables will be interpreted as a morphism $M^n \to M$, while a formula with $n$ free variables will be interpreted as a morphism $M^n \to \SC$. We assume the variables are numbered $x_1, x_2, \ldots$; let $v(t)$ list the free variables appearing in a term $t$, and $v(\varphi)$ those appearing in a formula $\varphi$.

\begin{definition}
	Each $\Ll$-term $t$ receives an interpretation $t^\M\colon M^n \to M$ in $\M$ as follows. Note that all constant symbols $c \in \Ll$ have already been given interpretations $c^\M\colon \tl \to \M$.
	\begin{enumerate}
		\item Each variable $x$ has interpretation $x^\M = \id_M\colon M \to M$.
		
		\item For $\Ll$-terms $t_1,\ldots,t_n$ with interpretations $t_i^\M\colon M^{k_i} \to M$, the \textit{product term} $(t_1,\ldots,t_n)$ has interpretation $(t_1,\ldots,t_n)^\M$ the composite
		\begin{center}
		\begin{tikzcd}
			M^d \arrow[r,"\Delta_\sigma"] & M^k \arrow[rr,"t_1 \times \cdots \times t_k"] & & M^n
		\end{tikzcd}
		\end{center}
		where $k = \sum k_i$, and $\Delta_\sigma$ is the diagonal of signature $\sigma = (m_1,\ldots,m_k)$, where $(x_{m_1},\ldots,x_{m_k}) = \big( v(t_1),\ldots,v(t_n) \big)$.

		\item A function term $f(t_1,\ldots,t_n)$ has interpretation $[f(t_1,\ldots,t_n)]^\M$ the composite
		\begin{center}
		\begin{tikzcd}
			M^d \arrow[rr,"{(t_1,\ldots,t_n)^\M}"] & & M^n \arrow[r,"f^\M"] & M
		\end{tikzcd}
		\end{center}
	\end{enumerate}
\end{definition}
The point of $\Delta_\sigma$ is to ensure that variables shared between different terms are identified with each other. Before giving interpretations of $\Ll$-formulae, we need to define \textit{Boolean connectives} in a topos, which are maps $\SC^n \to \SC$, and \textit{quantifiers}, which are maps $\SC^M \to \SC$.

\begin{definition}[Boolean connectives {\cite[139]{goldblattTopoiCategorialAnalysis1984}}]\label{defn:bool-ctves}\
	\begin{enumerate}
		\item $\false\colon \tl \to \SC$ (called \textit{false}) is the character of $\unq\colon \intl \to \tl$, and $\lnot\colon \SC \to \SC$ is the character of $\false$. In sum, $\false$ and $\lnot$ are the unique maps making both squares pullbacks:
		\begin{center}
		\begin{tikzcd}
			\intl \arrow[r,"\unq"] \arrow[d,"\unq"'] & \tl \arrow[r,"\unq\, =\, \id_\tl"] \arrow[d,"\false"] & \tl \arrow[d,"\true"] \\
			\tl \arrow[r,"\true"] & \SC \arrow[r,"\lnot"] & \SC
		\end{tikzcd}
		\end{center}
	
		\item $\land\colon \SC^2 \to \SC$ is the character of $\true^2\colon \tl \to \SC^2$.
		
		\item By the definition of the coproduct $\SC\!+\!\SC$, there is a unique map $k\colon \SC\!+\!\SC \to \SC^2$ making the following commute:
		\begin{center}
		\begin{tikzcd}
			\SC \arrow[r,"i_1"] \arrow[dr,"\true \times \id_\SC"'] & \SC\!+\!\SC \arrow[d,dashed,"k"] & \SC \arrow[l,"i_2"'] \arrow[dl,"\id_\SC\! \times \true"] \\
			& \SC^2
		\end{tikzcd}
		\end{center}
		
		Then, $\lor\colon \SC^2 \to \SC$ is the character of $\im k$.\footnote{As in abelian categories, any arrow $f$ in a topos has a unique epi-monic factorisation $f = \im f\, \circ\ \coim f$.}
		
		\item $\to\colon \SC^2 \to \SC$ is the character of the equaliser of \begin{tikzcd}[cramped]
			\SC^2 \arrow[r,shift left,"\pi_1"] \arrow[r,shift right,"\land"'] & \SC
		\end{tikzcd}.
		
		\item $\leftrightarrow\colon \SC^2 \to \SC$ is the character of the diagonal $\Delta\colon \SC \to \SC^2$.
	\end{enumerate}
\end{definition}

Definition \ref{defn:bool-ctves} is the appropriate generalisation of the usual Boolean connectives in $\Set$.

\begin{definition}[Quantifiers {\cite[245]{goldblattTopoiCategorialAnalysis1984}}]\label{defn:quant}
	Fix an object $M \in \C$.
	\begin{enumerate}
		\item By Definition \ref{defn:expon}, the composite \begin{tikzcd}
			M \arrow[r,"\unq"] & \tl \arrow[r,"\true"] & \SC
		\end{tikzcd} corresponds to a map $\bar{a}\colon \tl \to \SC^M$, implicitly using the natural isomorphism $\tl \!\times\! M \cong M$. Then, $\forall_M\colon \SC^M \to \SC$ is the character of $\bar{a}$.
		
		\item Let $\ni_M\colon K \mon \SC^M \!\times\! M$ be as in Definition \ref{defn:power-obj}. Then, $\exists_M\colon \SC^M \to \SC$ is the character of $\pi_1\, \circ \ni_M$.
	\end{enumerate}
\end{definition}

In $\Set$, $\forall_M, \exists_M\colon \Pow(M) \to \textbf{2}$ are such that $\forall_M(B) = 1 \iff B = M$ and $\exists_M(B) = 1 \iff B \neq \varnothing$. Definition \ref{defn:quant} generalises this to any topos. Now we can define interpretations $\varphi^\M$ of formulae $\varphi$:

\begin{definition}
	Every $\Ll$-formula $\varphi$ has an interpretation $\varphi^\M \colon M^n \to \SC$ as follows:
	\begin{enumerate}
		\item For $\Ll$-terms $s$, $t$, the formula $(s=t)$ has interpretation $(s=t)^\M$ the composite
		\begin{center}
		\begin{tikzcd}
			M^d \arrow[r,"{(s,t)^\M}"] & M^2 \arrow[r,"\chi_\Delta"] & \SC
		\end{tikzcd}
		\end{center}
		where $\chi_\Delta$ is the character of the usual diagonal $\Delta\colon M \to M^2$.
	
		\item For $\Ll$-terms $t_1, \ldots, t_n$ and a relation symbol $R \in \Ll$, the formula $R(t_1,\cdots,t_n)$ has interpretation $[R(t_1,\ldots,t_n)]^\M$ the composite
		\begin{center}
		\begin{tikzcd}
			M^d \arrow[rr,"{(t_1,\ldots,t_n)^\M}"] & & M^n \arrow[r,"R^\M"] & M
		\end{tikzcd}
		\end{center}
	
		\item For an $\Ll$-formula $\varphi$, the formula $\lnot \varphi$ has interpretation $(\lnot \varphi)^\M$ the composite
		\begin{center}
		\begin{tikzcd}
			M^n \arrow[r,"\varphi^\M"] & \SC \arrow[r,"\lnot"] & \SC
		\end{tikzcd}
		\end{center}
	
		\item For $\Ll$-formulae $\varphi$, $\psi$ with interpretations $\varphi^\M\colon M^n \to \SC$, $\psi^\M\colon M^k \to \SC$, the formula $\varphi \land \psi$ has interpretation $(\varphi \land \psi)^\M$ the composite
		\begin{center}
		\begin{tikzcd}
			M^d \arrow[r,"\Delta_\sigma"] & M^{n+k} \arrow[rr,"\varphi^\M \times \psi^\M"] & & \SC^2 \arrow[r,"\land"] & \SC
		\end{tikzcd}
		\end{center}
		where $\Delta_\sigma$ has signature $\sigma = (m_1,\ldots,m_{n+k})$ for $(x_{m_1},\ldots,x_{m_{n+k}}) = \big( v(\varphi),v(\psi) \big)$. $(\varphi \lor \psi)$, $(\varphi \to \psi)$, $(\varphi \liff \psi)$ are interpreted similarly.
		
		\item For an $\Ll$-formula $\varphi$ with interpretation $\varphi^\M\colon M^{n+1} \to \SC$, the formula $\forall x\, \varphi$ is interpreted as follows. As in Definition \ref{defn:expon}, $\varphi^\M$ corresponds to a map $\bar{\varphi}_x\colon M^n \to \SC^M$, where we ``curry'' in the position corresponding to the variable $x$. Then, $(\forall x\, \varphi)^\M$ is the composite
		\begin{center}
		\begin{tikzcd}
			M^n \arrow[r,"\bar{\varphi}_x"] & \SC^M \arrow[r,"\forall_M"] & \SC
		\end{tikzcd}
		\end{center}
	
		\item Similarly, $\exists x\, \varphi$ has intepretation $(\exists x\, \varphi)^\M$ the composite
		\begin{center}
		\begin{tikzcd}
			M^n \arrow[r,"\bar{\varphi}_x"] & \SC^M \arrow[r,"\exists_M"] & \SC
		\end{tikzcd}
		\end{center}
	\end{enumerate}
\end{definition}

\begin{definition}\label{defn:models}
	Let $\varphi$ be an $\Ll$-formula with interpretation $\varphi^\M\colon M^n \to \SC$, and $a_1, \ldots, a_n\colon \tl \to M$. We say $\M \models \varphi(a_1,\ldots,a_n)$ if the following diagram commutes:
	\begin{center}
	\begin{tikzcd}[column sep=tiny]
		& M^n \arrow[rd,"\varphi^\M"] & \\
		\tl \arrow[ru,"a_1 \times \cdots \times a_n"] \arrow[rr,"\true"'] & & \SC
	\end{tikzcd}
	\end{center}

In particular, for an $\Ll$-sentence $\varphi$, $\M \models \varphi \iff \varphi^\M = \true$.
\end{definition}

The usual internal logic of a topos is \textit{intuitionistic logic}, where the law of excluded middle $\varphi \lor \lnot \varphi \equiv \top$ may fail. Definition \ref{defn:models} gives some insight into why this might be the case. Just because $\varphi^\M \neq \true$ doesn't necessarily mean $(\lnot\varphi)^\M = \true$; $(\lnot\varphi)^\M$ could be equal to some other ``truth value'' $\tl \to \SC$.

\begin{definition}
	The definition of the coproduct $\tl\!+\!\tl$ gives a unique map $k\colon \tl\!+\!\tl \to \SC$ making the following commute:\vspace*{-2mm}
	\begin{center}
	\begin{tikzcd}
		\tl \arrow[r,"i_1"] \arrow[dr,"\true"'] & \tl\!+\!\tl \arrow[d,dashed,"k"] & \tl \arrow[l,"i_2"'] \arrow[dl,"\false"] \\
	& \SC
	\end{tikzcd}
	\end{center}

	A topos $\C$ is \textit{Boolean} or \textit{classical} if $k$ is an isomorphism.
\end{definition}

So a Boolean topos is one where the subobject classifier $\SC$ is essentially the coproduct $\true\!+\!\false$. These are exactly the topoi in which the law of excluded middle always holds, i.e.\ the internal logic is classical. However, it must be emphasised that this is not the general case.

\begin{example}
	The \textit{arrow category} $\Set^\to$ is that whose objects are functions $f\colon A \to B$ (in $\Set$), and morphisms from $f\colon A \to B$ to $g\colon C \to D$ are commutative squares:\vspace*{-2mm}
	\begin{center}
	\begin{tikzcd}
		A \arrow[d,"f"'] \arrow[r] & C \arrow[d,"g"] \\
		B \arrow[r] & D
	\end{tikzcd}
	\end{center}
	$\Set^\to$ forms a topos whose internal logic is \textit{three-valued}. The three truth values have a natural interpretation as a time-like logic: fixing a point in time $t_0$, the truth values are \textit{always true} ($\true$), \textit{always false} ($\false$), and \textit{false before $t_0$, true afterwards} ($\ch$). The truth tables for the Boolean connectives are as follows (the row is the first argument).\\[-2mm]
	
	\hfill
	$\begin{array}{c|ccc}
		\land & \true & \ch & \false \\
		\hline
		\true & \true & \ch & \false \\
		\ch & \ch & \ch & \false \\
		\false & \false & \false & \false \\
	\end{array}$
	\hfill
	$\begin{array}{c|ccc}
		\lor & \true & \ch & \false \\
		\hline
		\true & \true & \true & \true \\
		\ch & \true & \ch & \ch \\
		\false & \true & \ch & \false \\
	\end{array}$
	\hfill
	$\begin{array}{c|ccc}
		\to & \true & \ch & \false \\
		\hline
		\true & \true & \ch & \false \\
		\ch & \true & \true & \false \\
		\false & \true & \true & \true \\
	\end{array}$
	\hfill\
	
\end{example}

%
%
%
%
%
%
%
%
%
%
%

\section{Categorical set theory}

Having formalised model theory within a topos, we now turn to set theory. The reader may have previously heard the vague assertion that ``category theory can serve as an alternative foundation for mathematics''. Here, we give one way to make this precise.

There is a key difference between \textit{material} (classical) and \textit{structural} (categorical) set theory. In material set theory, sets are defined by their elements; for example, the sets $(A\!\times\!A)\!\times\!A$ and $A\!\times\!(A\!\times\!A)$, although (naturally) isomorphic, are different since they have different elements. Structural set theory is ``set theory up to isomorphism''---any two objects with an isomorphism between them should be considered ``the same''.

In this section, we will examine which axioms of $\ZFC$ are true in a topos.
Indeed, Theorem \ref{thm:topos-power-obj} says exactly that any topos satisfies the power set axiom. 
Furthermore, \cite[\S4.8]{goldblattTopoiCategorialAnalysis1984} shows how the existence of a subobject classifier can be construed as a comprehension principle, hence any topos satisfies the axiom of comprehension.\footnote{The situation with regularity is complex; \cite{awodey14} shows that any topos can be construed both as a model of foundation, \textit{and} as a model of anti-foundation. \cite{bauer20} gives a good explanation.}

So what's missing? General topoi do not have to be extensional: two distinct morphisms could act the same with regards to composition.

\begin{definition}
	An \textit{element} of $X \in \C$ is a map $x\colon \tl \to X$.
\end{definition}

\begin{definition}
	A topos is \textit{well-pointed} if for every $f,g\colon A \to B$, whenever $fx = gx$ for all elements $x\colon \tl \to A$, then $f=g$.
\end{definition}

So, a well-pointed topos is one satisfying extensionality. We still don't have the axiom of infinity---the category of finite sets is a well-pointed topos. To amend this, we introduce an ``infinite'' object which behaves like the natural numbers.

%
%

\begin{definition}\label{defn:nno}
	A \textit{natural numbers object} (NNO) is a diagram
	\begin{tikzcd}
	\tl \arrow[r,"z"] & \N \arrow[r,"s"] & \N
	\end{tikzcd}
	such that for every diagram
	\begin{tikzcd}
	\tl \arrow[r,"c"] & A \arrow[r,"v"] & A
	\end{tikzcd}
	there is a unique arrow $h\colon \N \to A$ making the following commute:\vspace*{-4mm}
	\begin{center}
	\begin{tikzcd}[row sep=small]
		& \N \arrow[r,"s"] \arrow[dd,"h"] & \N \arrow[dd,"h"] \\
		\tl \arrow[ru,"z"] \arrow[rd,"c"'] \\
		& A \arrow[r,"v"] & A
	\end{tikzcd}
	\end{center}
\end{definition}

Intuitively, we should think of $\N$ as analogous to the natural numbers, $z\colon \tl \to \N$ as analogous to the number zero, and $s\colon \N \to \N$ as analogous to the successor function. Then, Definition \ref{defn:nno} says that, given $c \in A$ and $v\colon A \to A$, we can define a map $h\colon \N \to A$ ``by recursion'', i.e. $h(0)=c$, $h(n+1) = v(h(n))$, and this gives a unique function $h$.\footnote{Definition \ref{defn:nno} does not claim that \textit{all} functions $h\colon \N\to A$ are defined this way.}

What about the axiom of choice? The standard formulation is that for any collection $\{ A_i: i \in I \}$ of sets, there is a sequence $(a_i)_{i \in I}$ with each $a_i \in A_i$. But an indexed collection $\{ A_i: i \in I \}$ is simply a surjective function $f\colon A \to I$, where $A = \bigsqcup_{i \in I} A_i$, and a sequence $(a_i)_{i \in I}$ is simply a function $g\colon I \to A$. The requirement that $a_i \in A_i$ is equivalent to saying that $fg = \id_I$. This leads to a categorical formulation:

\begin{definition}
	A topos $\C$ satisfies the \textit{axiom of choice} if any epi $f\colon A \to I$ in $\C$ has a right inverse $g\colon I \to A$, i.e.\ $fg = \id_I$.
\end{definition}


At this stage, we have almost all of the $\ZFC$ axioms. The first attempt to axiomatise mathematics with category theory was Lawvere's $\ETCS$.
In modern language, $\ETCS$ says the category of sets is a well-pointed topos with an NNO, satisfying the axiom of choice. This is equiconsistent with $\BZC$ (bounded Zermelo with choice), which is $\ZFC$ without replacement, and where comprehension is restricted to formulae having bounded quantifiers.

To get the full strength of $\ZFC$, we need to add \textit{replacement}. \cite[\S 8]{mclartyExploringCategoricalStructuralism2004} shows one way to do this. The \textit{language of categories} has two sorts: objects $X$ and morphisms $f$, with symbols $\{ \operatorname{source}, \operatorname{range}, \circ \}$ of the appropriate signature. Suppose $R(f,X)$ is a predicate definable in this logic, and $A \in \C$ is such that for any $x\colon \tl \to A$, there is unique (up to isomorphism) $S_x \in \C$ with $R(x,S_x)$. Then, \textit{replacement} asserts that there is $g\colon S \to A$ such that every $S_x$ is obtained from $x$ by pullback along $g$:
\begin{center}
\begin{tikzcd}
	S_x \arrow[r,hook] \arrow[d,"\unq"] & S \arrow[d,"g"] \\
	\tl \arrow[r,"x"] & A
\end{tikzcd}
\end{center}

The replacement axiom $\R$ asserts that this is true for every such predicate $R(f,X)$ and object $A \in \C$. \cite{mclartyExploringCategoricalStructuralism2004} shows that $\ETCS\!+\!\R$ is equivalent to $\ZFC$.

We can also formalise additional set theory axioms (such as large cardinal axioms) in a topos. Since every set admits a bijection to a unique cardinal, from the point of view of structural set theory, a cardinal is just a set. For objects $X, Y \in \C$, we say $X \leq Y$ if there is a monic $X \to Y$, and say $X < Y$ if not $Y \leq X$. Then, we can give the usual definitions of uncountable sets, limit cardinals, etc. For example, here is one way to define inaccessible cardinals in a topos.

\begin{definition}\
	\begin{enumerate}
		
		\item Objects $X \in \C$ induce discrete categories $\bar{X} = \{ x: \tl \to X \}$. $X \in \C$ is \textit{regular} if for every $Y < X$ and diagram $\D\colon \bar{Y} \to \C$ with colimit $Z$, if $\D y < X$ for all $y \in \bar{Y}$, then $Z < X$.
		
		\item $X$ is \textit{inaccessible} if it is an uncountable, regular, strong limit.
	\end{enumerate}
\end{definition}



\printbibliography

\end{document}